\documentclass[a4paper]{article}
\usepackage[T1]{fontenc}
\usepackage[ansinew]{inputenc}
\usepackage{psfrag}
\usepackage{amsmath}
\usepackage{graphicx}
\usepackage[colorlinks]{hyperref}

\begin{document}
\begin{flushright}
project\end{flushright}
\begin{center}
\title*{\Large{\textbf{
A hierarchical technique for estimating location parameter in the
presence of missing data}\\}}
\end{center}
\begin{center}
\author*{\large{\textbf{By SERGEY TARIMA}}} \\
\small{Department of Statistics, University of Kentucky, \\
Lexington, Kentucky, 40506-0027, U.S.A. \\ stari@ms.uky.edu}
\end{center}

\begin{center}
\author*{\large{\textbf{YURIY DMITRIEV}}} \\
\small{Department of Applied Mathematics and Cybernetics, Tomsk
State University, \\
Tomsk, 634050, Russia \\ dmit@fpmk.tsu.ru }
\end{center}

\begin{center}
\author*{\large{\textbf{RICHARD KRYSCIO}}} \\
\small{Department of Statistics and Department of Biostatistics, University of Kentucky, \\
Lexington, Kentucky, 40506-0027, U.S.A. \\ kryscio@email.uky.edu}
\end{center}

\section*{\begin{center} Summary \end{center}}

This paper proposes a hierarchical method for estimating the
location parameters of a multivariate vector in the presence of
missing data. At $i^{th}$ step of this procedure an estimate of the
location parameters for non-missing components of the vector is
based on combining the information in the subset of observations
with the non-missing components with updated estimates of the
location parameters from all subsets with even more missing
components in an iterative fashion. If the variance-covariance
matrix is known, then the resulting estimator is unbiased with the
smallest variance provided missing data are ignorable. It is also
shown that the resulting estimator based on consistent estimators of
variance-covariance matrices obtains unbiasedness and the smallest
variance asymptotically. This approach can also be extended to some
cases of non-ignorable missing data. Applying the methodology to a
data with random dropouts yields the well known Kaplan-Meier
estimator.\\

{Some key words: Parameter estimation; Missing data; Hierarchical
technique for missing data}

\section*{\begin{center} 1. Introduction \end{center}}

Censored and missing data are unavoidable parts of many rectangular
data sets. For the purposes of handling these kind of data many
different approaches have been developed in recent years. Little and
Rubin (2002) considered a taxonomy of missing-data methods
consisting of procedures based on completely recorded units,
weighting procedures, imputation-based and model-based procedures.
All these procedures can be classified into two general categories:
imputational and non-imputational techniques.

The first category contains a variety of single and multiple
imputation methods including mean substitution, last observation
carried forward, and imputational techniques for likelihood-based
approaches. Multiple Imputation (MI) (Rubin, 1987) is now the
accepted standard with several statistical packages supplying easy
to use software for applying this method (see, for example,
procedures MI and MIANALYZE in SAS, 2002). Monte Carlo Markov Chain
(MCMC) provides a flexible tool for MI. Some illustrative MCMC
examples are described by Schafer (1997).  Expectation Maximization
algorithm (Dempster, Laird, Rubin, 1977) for maximum likelihood
estimators and approximate Bayesian Bootstrap (Rubin and Schenker,
1986) for stratified samples are in this category. In addition,
several authors have investigated the small sample as well as large
sample properties of estimators based on multiple imputation
(Barnard and Rubin, 1999).

The second category consists of non-imputational techniques with the
complete case method and available case method being the most
popular (Verbeke and Molenberghs, 2000). In addition considerable
methodology has been constructed for obtaining maximum likelihood
estimators: parameter estimation on incomplete data in general
linear models (Ibrahim, 1990); pattern set mixture models (Little,
1993), including the analysis based on pattern mixture models and
selection models. The analysis based on pattern mixture models is
the one in which inference for a function of the location parameters
is obtained by combining in some weighted fashion estimates obtained
from each pattern of missing components observed in the data
(Molenberghs, Michiels, Kenward, Diggle, 1998). Pattern mixture
models are the closest analogues to the technique proposed in this
paper, but the proposed method does not depend on assuming a
parametric family.

To develop a new distribution free non-imputation approach for
estimation on missing data we reviewed some methods proposed and
developed for involving auxiliary information in statistical
function estimation. One important method due to Pugachev (1973) is
the method of correlated processes which uses correlation effect
between auxiliary information and empirical data for incorporating
auxiliary information in statistical estimation. This method was
later developed and extended by Gal'chenko and Gurevich (1991) who
incorporated the estimators from previous experiments into the
current estimator. The estimators obtained by these approaches
provide smaller or asymptotically smaller variances than the
variance of the current estimator. The further extension which is
the subject of this paper provides a methodological basis for
statistical estimation for missing data.

This new method is introduced in Section 2 and the asymptotic
properties of this method are then derived in Appendix. Applications
to the situation where missing data is due to right censoring is
considered in Section 3 and shows that in this important special
case the method produces the well known Kaplan-Meier estimator. The
other applications to samples from a bivariate random variable with
ignorable and a special case of non-ignorable missing data are
presented in Section 4. In this section considered the vector of
means estimation at general pattern of ignorable missing data and
change score estimation at random dropout. Conclusions are stated in
Section 5.

\section*{\begin{center} 2. Methodology \end{center}}

 \subsection*{2.1. Notation}

Suppose $\mathbf{X}_1, ... , \mathbf{X}_N$ are independent and
identically distributed random vectors with common probability
distribution $P_{\mathbf{X}}(\mathbf{x})$, where $\mathbf{x} \in
\mathcal{X} \subset \mathcal{R}^{Q}$, $N$ and $Q$ are finite and
strictly positive integers. But $\mathbf{X}_1,... , \mathbf{X}_N$
are not observed directly. These data are subject to a missing data
mechanism by corresponding vectors indicating nonresponse:
$\mathbf{R}_1,...,\mathbf{R}_N$. Here
$\mathbf{R}_n=\left(R_{n1},...,R_{nQ}\right)$ and $R_{nq} \in
\{0,1\}$, $n=1,...,N$, $q=1,...,Q$. In the notation $R_{nq}=1$
indicates response and $R_{nq}=0$ indicates non-response. What is
really observed is a random vector $\mathbf{Y}_1, ... ,
\mathbf{Y}_N$, where $\mathbf{Y}_n=\left(Y_{n1},...,Y_{nQ}\right)$,
$Y_{nq}=X_{nq}$ if $R_{nq} = 0$ and $Y_{nq}$ is missing if $R_{nq} =
1$, $n=1,...,N$, $q=1,...,Q$.

Let $\Theta=\left(\theta_{1},...,\theta_{S}\right)$ take values in
$\mathcal{R}^S$, where $\theta_{s} =
\int_{\mathcal{X}}\varphi_{s}(\mathbf{x})
dP_{\mathbf{X}}(\mathbf{x})$ with $\varphi_{s}(\mathbf{x})$ a known
function defined on $\mathcal{R}^Q$, $\theta_{s} \in \mathcal{R}$,
$s=1,...,S$.

Several examples of
$\varphi_{s}(\mathbf{x})=\varphi_{s}({x_1,...,x_Q})$ follow.

In this paper the location parameter estimation is emphasized. If
$\varphi_{s}({x_1,...,x_Q}) = x_1$, then $\theta_{s}$ is a mean of
$x_1$. If $\varphi_{s}({x_1,...,x_Q})$ is an indicator function of
some event defined by the variables $x_1,...,x_Q$, then the
parameter $\theta_{s}$ becomes the probability of this event. Hence,
a Cumulative Distribution Function ($CDF$) can be estimated. The
obtained $CDF$ estimator can further be used to estimate
percentiles, median, interquartile range, and many other parameters.

This approach is not restricted only to location parameter
estimation. If $\varphi_{s}({x_1,...,x_Q}) = x_1 x_Q$ then
$\theta_{s}$ is a mixed moment of $x_1$ and $x_Q$. In general, we
are not excluding from consideration the possibility of more
intricate forms for $\varphi_{s}({x_1,...,x_Q})$, for example
$\varphi_{s}({x_1,...,x_Q})= x_Q \log \left({x_1 x_2}\right)$.

Though the location parameter estimation is the main objective of
this paper, the methodology presented in this section accommodates
all these cases.

First, consider an ignorable mechanism of missing data generation.
The idea of how to apply this approach to non-ignorable missing data
is considered in subsection 2.4 with a special case in Section 4.

 \subsection*{2.2. Hierarchical structure}

Let $R_{ij}$ denote an indicator vector having exactly $(i-1)$ zeros
for $i=1,...,Q$. For a given $i$ we have $j=1,...,\left(_i^Q\right)$
different patterns with exactly $(i-1)$ zeros. Let $J_{ij}$ denote
the subsample size for the $i^{th}$ level and the $j^{th}$ pattern,
where $J_{ij} \ge 0 $. Let $\Theta_{ij}$ denote the subset of the
$S$ parameters $\theta_1,...,\theta_S$ which is estimable using only
the observations having the missing pattern defined by $R_{ij}$. Let
$\hat\Theta_{ij}$ denote this sample estimate assuming that $J_{ij}
\ge 0$. Notice that the $R's$ and corresponding estimates can be
arranged into a hierarchical structure as $i$ increases.

Example. If $Q=3$, then this hierarchical structure follows.

\begin{itemize}
    \item the subsample which contains complete observations defines
    the first level or root level $(i=1)$ and corresponds to the
    indicator vector (1,1,1);
    \item up to three subsamples define the second level $(i=2)$ and
    correspond to the three missing patterns $(1,1,0)$, $(1,0,1)$, and
    $(0,1,1)$;
    \item up to three subsamples define the third level $(i=3)$ and
    correspond to the three missing patterns $(1,0,0)$, $(0,1,0)$,
    and $(0,0,1)$.
\end{itemize}

We now use this hierarchy to improve the estimator $\hat\Theta_{ij}$
by using the information about the unknown value of $\Theta_{ij}$
from the next higher level. The improved estimator is

\begin{equation}
\tilde \Theta_{ij} = \hat \Theta_{ij} - \mathbf{K}_{ij} \left(
\mathbf{K}_{ij}^{*} \right)^{-1} \left( \hat{\cal{B}}_{ij} -
\tilde{\cal{B}}_{ij} \right).
\end{equation}

The elements of the $\mathbf{K}$ matrices and $\cal{B}$ vectors in
equation (1) are defined below. Assume there are $S^{*} \equiv
S(i,j)$ elements in $\Theta_{ij}$ and without loss of generality
assume these are numbered $1,...,S^{*}$. That is, assume,
$\Theta_{ij} = \left( \theta_1,...,\theta_{S^{*}} \right)$. Then
$\hat{\cal{B}}_{ij}=\left(
\hat{\cal{B}}_{ij1},...,\hat{\cal{B}}_{ijS^{*}} \right)$ and
$\tilde{\cal{B}}_{ij}=\left(
\tilde{\cal{B}}_{ij1},...,\tilde{\cal{B}}_{ijS^{*}} \right)$. To
define these vectors let ${\cal{B}}_{ijk}$ represent the subvector
of $\Theta_{ij}$ with its $k^{th}$ component missing for
$k=1,...,S^{*}$. Two estimates of ${\cal{B}}_{ijk}$ are computed
from the data. The first is based on the subsample defined by
$R_{ij}$; this is $\hat{\cal{B}}_{ijk}$ which is the subvector of
$\hat\Theta_{ij}$ with its $k^{th}$ component missing. The second is
based on data collected at the $(i+1)^{st}$ level, i.e.
${\tilde{\cal{B}}}_{ijk}$. It is possible that there are no
observations in the latter subsample in which case the corresponding
subvector is dropped from both ${\hat{\cal{B}}}_{ij}$ and
${\tilde{\cal{B}}}_{ij}$. The rectangular matrix $\mathbf{K}_{ij}$
is a block matrix defined as follows:

$$
\mathbf{K}_{ij} = \left \| Cov \left( \hat\Theta_{ij},
\hat{\cal{B}}_{ijk} \right) \right \|_{k=1,...,S^{*}}.
$$

\noindent The square block matrix $\mathbf{K}^{*}_{ij}$ is defined
as follows:

$$
\mathbf{K}^{*}_{ij} = \left \| Cov \left( \hat{\cal{B}}_{ijl},
\hat{\cal{B}}_{ijk} \right) + I_{\left[l=k\right]} Cov \left(
\tilde{\cal{B}}_{ijl}, \tilde{\cal{B}}_{ijk} \right) \right
\|_{k,l=1,...,S^{*}}.
$$

\noindent The estimator (1) defines the estimator with a
variance-covariance matrix

\begin{equation}
Cov \left( \tilde \Theta_{ij},  \tilde \Theta_{ij} \right)  = Cov
\left ( \hat \Theta_{ij}, \hat \Theta_{ij} \right)  -
\mathbf{K}_{ij} \left( \mathbf{K}_{ij}^{*} \right)^{-1}
\mathbf{K}^{T}_{ij}
\end{equation}

\noindent defining the smallest dispersion ellipsoid in a class

\begin{equation}
\tilde \Theta^{\Lambda}_{ij} = \hat \Theta_{ij} - \Lambda_{ij}
\left( \hat{\cal{B}}_{ij} - \tilde{\cal{B}}_{ij} \right)
\end{equation}

\noindent with respect to different choices of the matrix
$\Lambda_{ij}$ of proper dimensions. The estimators $\tilde
\Theta^{\Lambda}_{ij}$ define a class of unbiased estimators of
$\Theta_{ij}$.

In practice the true values of $\mathbf{K}_{ij}$,
$\mathbf{K}^{*}_{ij}$ and $Cov \left ( \hat \Theta_{ij}, \hat
\Theta_{ij} \right)$ usually are not available, in which case their
consistent estimators $\hat{\mathbf{K}}_{ij}$,
$\hat{\mathbf{K}}^{*}_{ij}$, and $\widehat{Cov} \left ( \hat
\Theta_{ij}, \hat \Theta_{ij} \right)$ are used instead.

This substitution modifies (1) and (2) to the following equations

\begin{equation}
\hat{\tilde{\Theta}}_{ij} = \hat \Theta_{ij} - \hat{\mathbf{K}}_{ij}
\left( \hat{\mathbf{K}}_{ij}^{*} \right)^{-1} \left(
\hat{\cal{B}}_{ij} - \tilde{\cal{B}}_{ij} \right)
\end{equation}

\noindent with

\begin{equation}
\widehat{Cov} \left( \tilde\Theta_{ij},  \tilde\Theta_{ij} \right) =
\widehat{Cov} \left ( \hat\Theta_{ij}, \hat\Theta_{ij} \right) -
\hat{\mathbf{K}}_{ij} \left( \hat{\mathbf{K}}_{ij}^{*} \right)^{-1}
\hat{\mathbf{K}}^{T}_{ij}.
\end{equation}

In addition to $det \left( \mathbf{K}^{*}_{ij} \right) > 0$ a new
requirement comes from (4) and (5): $\mathbf{\hat{K}}^{*}_{ij}$
should be positive definite. From $det \left( \mathbf{K}^{*}_{ij}
\right) > 0$ conclude that there exists a sufficiently large sample
size $N$ such that for any $n > N$ have $det
\left(\mathbf{\hat{K}}^{*}_{ij}\right) > 0$ with probability one.

\subsection*{2.3. Assumptions}

In order to obtain the unbiased estimator defined by (1) with the
smallest dispersion ellipsoid defined by (2) we need (for every
$ij$-subsample):

\begin{itemize}
    \item to know $K^{*}_{ij}$ and it should be positive definite,
    \item to know $K_{ij}$ (in many cases $K_{ij}$ consists of the elements of
$K^{*}_{ij}$ ),
    \item $E \left( \hat\Theta_{ij} \right) = \Theta_{ij}$, and
    \item $E \left( \hat{\cal{B}}_{ij} \right) = E \left( \tilde{\cal{B}}_{ij}
\right) = {\cal{B}}_{ij}$.
\end{itemize}

When $K^{*}_{ij}$ and $K_{ij}$ are not known their consistent
estimators provide unbiasedness and the smallest dispersion
asymptotically.

According to Little and Rubin (2002, p. 119) a missing-data
mechanism is ignorable if (1) the missing data are missing at random
and (2) parameters managing $X$ and $R$ are distinct that is in
different parameter spaces.

In case of ignorable missing data the missing data mechanism
splitting the original sample into subsamples is independent from
vector $\Theta$.

Hence, the methodology proposed in Section 2.2 can be applied to
ignorable missing data.

\subsection*{2.4. Adjustment for non-ignorable missing data
mechanism}

What does happen when missing-data mechanism is not ignorable? In
this case it is reasonable to assume that some or all of $S^{*}$
components of the vector $E \left( \hat{\cal{B}}_{ij} \right)$
differ from these of $E \left( \tilde{\cal{B}}_{ij} \right)$. In the
other words the bias was brought by missing data.

Suppose that missing data mechanism is managed not only by
parameters which are distinct from $\Theta_{ij}$ but also by $W$
parameters which are defined in a parameter space of $\Theta_{ij}$ .
If $W < S^{*}$, then it can be expected that there exist $S^{*}-W$
parameters independent from the missing data mechanism and they can
be used as a components of the vectors $\hat{\cal{B}}_{ij}$ and
$\tilde{\cal{B}}_{ij}$.

Hence, the purpose is \textbf{to find the parameters independent
from the missing data mechanism}. And use these in formulas
(1),(2),(4) and (5). Example of such a case is considered in Section
4.

In order to illustrate applicability of the methodology described in
the section consider the following special case.

\section*{\begin{center} 3. Random Dropout \end{center}}

Right censored data is one of the most common problems statisticians
face. This problem can be formulated in terms of missing data with
monotone missing data structure.

Suppose $X_1, \dots, X_N$ are independent and identically
distributed random variables with an unknown cumulative distribution
function $F(t)$, $t \in [0,\infty)$. But $X_1, \dots, X_N$ are not
observed directly since some are distorted by $M_1,...,M_N$
generated by a random missing mechanism. The observed sample is
$Y_1, \dots, Y_N$, where $Y_n = X_n$ if $M_n \geq X_n$ and $Y_n =
M_n$ otherwise, $n \in \{1,...,N\}$.

Assume the observed events occur at $t_1 < ... < t_S$, where $S \leq
N$. Consider an arbitrary event time $t_s$. On the basis of complete
(not censored at or before $t_s$) observations the empirical
estimators $\hat{F}\left(t_s\right)$ and
$\hat{F}\left(t_{s-1}\right)$ can be calculated. In addition to
these estimators an estimator $\tilde{F}\left(t_{s-1}\right)$ was
obtained on the basis of the data independent from complete
observations. At $s=2$ the estimator $\tilde{F}\left(t_{1}\right)$
uses only the observations censored at $t_2$. But at an arbitrary
$s^{th}$ step $\tilde{F}\left(t_{s-1}\right)$ represents an
estimator absorbing information from all previously censored
observations. We will not need to define its form explicitly because
in a recursive approach considered below we use the estimator
$\tilde{F}^{0}\left(t_{s-1}\right)$ absorbing information from
$\hat{F}\left(t_{s-1}\right)$ and $\tilde{F}\left(t_{s-1}\right)$.

From we (1) obtain the following equation

\begin{equation}
\tilde{F}^{0}(t_s) = \hat{F}(t_s) - \frac{Cov \left( \hat{F} \left(
t_s \right), \hat{F} \left( t_{s-1} \right) \right)} { Var \left(
\hat{F} \left( t_{s-1} \right) \right)+ Var \left( \tilde{F} \left(
t_{s-1} \right) \right)} \left[\hat{F}\left(t_{s-1}\right) -
\tilde{F}\left(t_{s-1}\right)\right].
\end{equation}

Considering the class of unbiased estimators
$$\tilde{F}^{\lambda}(t_{s-1}) = \hat{F}(t_{s-1}) - \lambda
\left[\hat{F}\left(t_{s-1}\right) -
\tilde{F}\left(t_{s-1}\right)\right],$$ the estimator

\begin{equation}
\tilde{F}^{0}(t_{s-1}) = \hat{F}(t_{s-1}) - \frac{Var \left( \hat{F}
\left( t_{s-1} \right) \right)} { Var \left( \hat{F} \left( t_{s-1}
\right) \right)+ Var \left( \tilde{F} \left( t_{s-1} \right)
\right)} \left[\hat{F}\left(t_{s-1}\right) -
\tilde{F}\left(t_{s-1}\right)\right]
\end{equation}

\noindent provides the smallest variance

\begin{equation}
Var \left(\tilde{F}^{0}(t_{s-1})\right) =
\frac{Var\left(\tilde{F}\left(t_{s-1}\right)\right)
Var\left(\hat{F}\left(t_{s-1}\right)\right)}
{Var\left(\hat{F}\left(t_{s-1}\right)\right)+
Var\left(\tilde{F}\left(t_{s-1}\right)\right)}.
\end{equation}

\noindent The estimator (7) can be rewritten as
$$\tilde{F}^{0}(t_{s-1}) = \hat{F}(t_{s-1})
\frac{Var\left(\tilde{F}\left(t_{s-1}\right)\right)} {Var
\left(\hat{F}\left(t_{s-1}\right)\right)+
Var\left(\tilde{F}\left(t_{s-1}\right)\right)} $$

\begin{equation}
+ \tilde{F}(t_{s-1})
\frac{Var\left(\hat{F}\left(t_{s-1}\right)\right)}
{Var\left(\hat{F}\left(t_{s-1}\right)\right)+
Var\left(\tilde{F}\left(t_{s-1}\right)\right)}.
\end{equation}

From (8) we have

\begin{equation}
Var\left(\tilde{F}(t_{s-1})\right) =
\frac{Var\left(\tilde{F}^{0}(t_{s-1})\right)
Var\left(\hat{F}\left(t_{s-1}\right)\right)}
{Var\left(\hat{F}\left(t_{s-1}\right)\right) -
Var\left(\tilde{F}^{0}(t_{s-1})\right) }.
\end{equation}

It is interesting to see that from (10) we can write
$$\left(Var\left(\tilde{F}^{0}(t_{s-1})\right)\right)^{-1} =
\left(Var\left(\hat{F}(t_{s-1})\right)\right)^{-1} +
\left(Var\left(\tilde{F}(t_{s-1})\right)\right)^{-1},$$ which shows
that Fisher information in $\tilde{F}^{0}(t_{s-1})$ is a sum of the
Fisher information in $\hat{F}(t_{s-1})$ and in
$\tilde{F}(t_{s-1})$.

Substituting (10) into (9) we obtain

$$\tilde{F}(t_{s-1}) =
\tilde{F}^{0}(t_{s-1})
\frac{Var\left(\hat{F}\left(t_{s-1}\right)\right)}
{Var\left(\hat{F}\left(t_{s-1}\right)\right)-
Var\left(\tilde{F}^{0}\left(t_{s-1}\right)\right)}$$

\begin{equation}
- \hat{F}(t_{s-1})
\frac{Var\left(\tilde{F}^{0}\left(t_{s-1}\right)\right)}
{Var\left(\hat{F}\left(t_{s-1}\right)\right)-
Var\left(\tilde{F}^{0}\left(t_{s-1}\right)\right)}.
\end{equation}

Applying the representation (11) of $\tilde{F}(t_{s-1})$ to the
equation $\tilde{F}^{0}(t_{s})$ have

\begin{equation}
\tilde{F}^{0}(t_{s})=\hat{F}(t_s)-
\frac{Cov\left(\hat{F}\left(t_{s}\right),\hat{F}\left(t_{s-1}\right)\right)}
{Var\left(\hat{F}\left(t_{s-1}\right)\right)}
\left[\hat{F}(t_{s-1})-\tilde{F}^{0}(t_{s-1})\right].
\end{equation}

\noindent Neither $\tilde{F}(t_{s-1})$ nor its variance appear in
(12) since the $\tilde{F}^{0}(t_{s-1})$ and its variance absorb all
information brought by $\tilde{F}(t_{s-1})$ and its variance.

Using the fact that $$Cov(\hat{F}(t_{s}),\hat{F}(t_{s-1}))=\frac
{F(t_{s-1})(1-F(t_{s}))} {n}$$ we have

\begin{equation}
\frac{Cov(\hat{F}(t_{s}),\hat{F}(t_{s-1}))}{Var(\hat{F}(t_{s-1}))} =
\frac{1-F(t_s)}{1-F(t_{s-1})}
\end{equation}

In $(13)$ the cumulative distribution function $F(\cdot)$ is not
known. Substituting its empirical estimator yields

$$ {\hat{\tilde{F}}}^{0}(t_s) = \hat{F}(t_s) -
\frac{1 - \hat{F}(t_s)} {1 - \hat{F}(t_{s-1})}
\left[\hat{F}(t_{s-1})-{\hat{\tilde{F}}}^{0}(t_{s-1})\right]
$$

\begin{equation}
= 1 - \frac{1 - \hat{F}(t_s)} {1 - \hat{F}(t_{s-1})}
\left[1-{\hat{\tilde{F}}}^{0}(t_{s-1})\right].
\end{equation}

From (14) have

\begin{equation}
1 - {\hat{\tilde{F}}}^{0}(t_s) =
\frac{1-\hat{F}(t_s)}{1-\hat{F}(t_{s-1})} \left[1 -
{\hat{\tilde{F}}}^{0}(t_{s-1})\right]
\end{equation}

The estimator ${\hat{\tilde{F}}}^{0}(t_{s-1})$ on the right side of
(15) was derived by applying $\hat{F}(\cdot)$ instead of unknown
$F(\cdot)$ (as it was done in (13)) on each of previous steps. Now
using survival function $S(\cdot)$ instead of $1-F(\cdot)$ the
equation (15) define the well-known Kaplan-Meier estimator (Kaplan
and Meier, 1958).

\section*{\begin{center} 4. Bivariate Case \end{center}}

Let $\mathbf{X}_1,...,\mathbf{X}_N$ be independent and identically
distributed random variables from a bivariate distribution with a
vector of means $\mu$ and a covariance matrix $\Sigma$, where
$\mathbf{X}_i=\left(X^{(1)}_i,X^{(2)}_i\right)$,
$\mu=\left(\mu_1,\mu_2\right)$, and
$\mathbf{\Sigma}=\left(%
\begin{array}{cc}
  \sigma_{11}^2 & \sigma_{12}^2 \\
  \sigma_{12}^2 & \sigma_{22}^2 \\
\end{array}%
\right)$ is a positive definite covariance matrix.

Applying the hierarchical structure developed in Section 2 we
summarize its content in the following table

\begin{tabular}{|c|c|c|c|c|c|}
  \hline
  Level $i$ & Subsample $j$ & $\mathbf{R}_{ij}$ & $J_{ij}$ & $\Theta_{ij}$ & Estimator \\
  \hline
  1 & 1 & (1,1) & $J_{11}$ & $\left( \mu_1, \mu_2 \right)^{T}$ & $\left( \bar{{X}}_{111},
  \bar{{X}}_{112}\right)^{T}$ \\
  2 & 1 & (1,0) & $J_{21}$ & $\mu_1$ & $\bar{{X}}_{211}$ \\
  2 & 2 & (0,1) & $J_{22}$ & $\mu_2$ & $\bar{{X}}_{222}$ \\
  \hline
\end{tabular}

The estimator of the vector $\left( \mu_1, \mu_2 \right)^{T}$ which
uses all information in the sample becomes
\begin{equation}
\left( \tilde{\mu}_1, \tilde{\mu}_2 \right)^{T} = \left(
\bar{{X}}_{111},\bar{{X}}_{112}\right)^{T} - \Lambda_0
\left(\bar{{X}}_{111}-\bar{{X}}_{211},
\bar{{X}}_{112}-\bar{{X}}_{222} \right)^{T}
\end{equation}

where \begin{equation}
 \Lambda_0 =
J_{11}^{-1}\left(\sigma_{11}^2,\sigma_{22}^2\right) \left(%
\begin{array}{cc}
  \sigma_{11}^2\left(1+\frac{J_{11}}{J_{21}}\right) &
  \sigma_{12}^2 \\
  \sigma_{12}^2 &
  \sigma_{22}^2\left(1+\frac{J_{11}}{J_{22}}\right) \\
\end{array}%
\right)^{-1}.\end{equation}

In a case when covariances in (17) are known the estimator (16) will
be unbiased with the smallest variance in class (3). If these
covariances are not known then their consistent estimates can be
used instead and the obtained estimator will not be the optimal
one anymore but it will converge to (16) in distribution
(see proposition 2 in Appendix). \\

An \textbf{ important special case } is $J_{22} = 0$. We discuss
this problem next.

\subsection*{ 4.1. Change Score Estimation}

Let $\delta=\mu_1-\mu_2$ be the change score we need to estimate.
This difference can be estimated with complete observations:
$\hat\delta = \bar{{X}}_{111}-\bar{{X}}_{112}$.

The estimator (1) takes the following form

\begin{equation}
\tilde{\delta} = \hat{\delta} - \frac{J_{11}}{J_{11}+J_{21}} \left(
1- \frac{\sigma^2_{12}}{\sigma^2_{11}}\right)
\left(\bar{{X}}_{111}-\bar{{X}}_{211} \right)
\end{equation}

\noindent with a variance

\begin{equation}
Var\left(\tilde{\delta}\right) = \frac{1}{J_{11}} \left(
\sigma^2_{11}- 2\sigma^2_{12}+\sigma^2_{22}  -
\frac{J_{21}}{\left(J_{11}+J_{21}\right)}\frac{\left(\sigma^2_{11}-\sigma^2_{12}\right)^2}{\sigma^2_{11}}\right).
\end{equation}

If $\sigma^2_{11} = \sigma^2_{12}$ then $\tilde{\delta} =
\hat{\delta}$ (the estimator based on \textbf{complete} cases).

If $\sigma^2_{12} = 0$ then $\tilde{\delta} =
\frac{J_{11}}{J_{11}+J_{21}} \bar{{X}}_{111} +
\frac{J_{21}}{J_{11}+J_{21}} \bar{{X}}_{211} - \bar{{X}}_{112}$ (the
estimator based on \textbf{available} cases).

\subsection*{ 4.2. Change Score Estimation at Compound Symmetry}

Let us assume $\mathbf{\Sigma}=\sigma^2\left(%
\begin{array}{cc}
  1 & \rho \\
  \rho & 1 \\
\end{array}%
\right)$ then

\begin{equation}
\tilde{\delta} = \hat{\delta} - \frac{J_{11}}{J_{11}+J_{21}} \left(
1- \rho\right) \left(\bar{{X}}_{111}-\bar{{X}}_{211} \right)
\end{equation}

\noindent with variance

\begin{equation}
Var\left(\tilde{\delta}\right) = \frac{\sigma^2}{J_{11}} \left(
2\left(1-\rho\right) - \frac{J_{21}}{J_{11}+J_{21}}\sigma^2 \left(1
- \rho\right)^2\right).
\end{equation}

\noindent If $\rho = 0$ then $\tilde{\delta} = \hat{\delta} -
\frac{J_{11}}{J_{11}+J_{21}} \left(\bar{{X}}_{111}-\bar{{X}}_{211}
\right)$ and $Var\left(\tilde{\delta}\right) =
\frac{\sigma^2}{J_{11}} \left( 2 -
\frac{J_{21}\sigma^2}{J_{11}+J_{21}} \right)$.

\noindent If $\rho = 1$ then $\tilde{\delta} = \hat{\delta}$ and
$Var\left(\tilde{\delta}\right) = 0$. \\

Now we return to the case where $J_{11}>0$, $J_{21}>0$, and
$J_{22}>0$ but assume data are not missing at random.

\subsection*{ 4.3. Non-ignorable Missing Data}

At \textbf{non-ignorable missing data} the parameters which do not
change after missing data transformations should be found. Let us
assume that the missing data case is the result of changed
experimental conditions, for example, $\Delta$ shift appears for
$X^{(1)}$ or $X^{(2)}$ if one of these components is missing. The
value of the $\Delta$ is unknown.

Using only incomplete observations obtain $\tilde{\delta} =
\bar{{X}}_{211}-\bar{{X}}_{222}$. In $\tilde{\delta}$ the $\Delta$
shift effect is canceled and $E\left(\tilde{\delta}\right) =
E\left(\hat{\delta}\right) = \Delta$. For these estimators
$Var\left(\hat\delta\right)=
J_{11}^{-1}\left(\sigma^2_{11}-2\sigma^2_{12}+\sigma^2_{22}\right)$,
$Var\left(\tilde\delta\right)=J_{22}^{-1}\sigma^2_{11}+J_{21}^{-1}\sigma^2_{22}$,
and the estimator (1) takes the following form
\begin{equation}
\tilde\delta^{\Lambda_0} = \hat\delta  -
\frac{\sigma^2_{11}-2\sigma^2_{12}+
\sigma^2_{22}}{\sigma^2_{11}\left(1+\frac{J_{11}}{J_{21}}\right)-2\sigma^2_{12}+
\sigma^2_{22}\left(1+\frac{J_{11}}{J_{22}}\right)} \left( \hat\delta
- \tilde\delta \right)
\end{equation}
with a variance
\begin{equation}
Var \left( \tilde\delta^{\Lambda_0} \right) =
\frac{\sigma^2_{11}-2\sigma^2_{12}+ \sigma^2_{22}}{J_{11}} -
\frac{\left(\sigma^2_{11}-2\sigma^2_{12}+
\sigma^2_{22}\right)^2}{J_{11}\left(\sigma^2_{11}\left(1+\frac{J_{11}}{J_{21}}\right)-2\sigma^2_{12}+
\sigma^2_{22}\left(1+\frac{J_{11}}{J_{22}}\right)\right)}.
\end{equation}

If the variances and covariances used in (22) and (23) are not
available, then their consistent estimators can be used. According
to Proposition 2 asymptotic properties continue to hold.

\section*{\begin{center} 5. Conclusion \end{center}}

If only the variance-covariance structure of a considered model is
known, the estimators proposed in this paper are unbiased and
provide the smallest variance in a class of unbiased estimators. In
the cases when one ought to estimate the parameters of variance
covariance structure with consistent estimators the estimators
obtain unbiasedness with the smallest variance asymptotically.

These estimators are not restricted to monotone missing data
structures and can be derived from the observations with a general
pattern of missing data. Despite the fact that these estimators are
obtained for the case of ignorable missing data they can also be
derived for some cases of non-ignorable mechanism of missing data. A
special case of nonignorable missing data considered in Section 5.

This approach does not require the assumptions on parametrical
families as many likelihood based methods and works when the first
two moments of the underlying distribution are finite.

Assuming asymptotical normality of the estimators obtained on
subsamples the final estimators obtained with proposed methodology
will be asymptotically normal as well. The two propositions in
Appendix provide asymptotical mean and variance for these
estimators.

Many standard statistical procedures may be used with these
estimators, for example, sample size determination or hypothesis
testing.

It was shown in section 3 that a well-known Kaplan-Meier estimator
is a result of applying this approach to right censoring data with
random dropout.

Overall, the nonparametric ground, the absence of any imputations in
any form, and the properties stated for finite and large sample
sizes make the proposed estimator distinct from the others and
applicable in many practical cases.

\section*{\begin{center} 6. References \end{center}}

\noindent Barnard, J., and Rubin, D.B. Small-sample degrees of
freedom with multiple imputation. {\it Biometrika}  86, 1999, no. 4, 948-955.\\

\noindent Casella, George; Berger, Roger L. Statistical inference.
{\it Duxbury}, 2002, 660pp.\\

\noindent Dempster, A. P., Laird, N. M., Rubin, D.B. Maximum
likelihood from incomplete data via the EM algorithm (with
discussion). {\it J. Roy. Statist. Soc.}
 B 39, 1977, 1-38.\\

\noindent Gal'chenko, M. V.; Gurevich, V. A. Minimum-contrast
estimation taking into account additional information. {\it
Journal of Soviet Math.}, 53, 1991, no. 6, 547-551.\\

\noindent Ibrahim G. Joseph. Incomplete data in generalized linear
models. {\it J. Amer. Statist. Ass.}, Vol. 85, No. 411, 1990,
765-769.\\

\noindent Kaplan E. L., P. Meier, Nonparametric estimator from
incomplete observations. {\it J. Amer. Statist. Ass.}, 53, 1958,
457-481.\\

Kulldorff, Gunnar. Contribution to the Theory of Estimation from
Grouped and Partially Grouped Samples. {\it Almqvist \& Wiksell},
Stockholm, 1961, pp. 142.

\noindent Little, R.J.A. Pattern-mixture models for multivariate
incomplete data. {\it J. Amer. Statist. Ass.}, 88, 1993,
125-134.\\

\noindent Little, R.J.A. and Rubin, D.B. Statistical Analysis with
missing data. {\it New York, Wiley-Interscience}, 2002\\

\noindent Molenberghs, G., Michiels, B., Kenward, M.G., and Diggle,
P.J. Missing data mechanisms and pattern-mixture models. {\it
Statistica Neerlandica}, 52, 1998, 153-161.\\

\noindent Pugachev, V. N. Mixed Methods of Determining Probabilistic
Characteristics. [in Russian] {\it Moscow, Soviet
Radio}, 1973, 256pp.\\

\noindent Rubin, D.B. Inference and missing data. {\it
Biometrika}, 63, 1976, no. 3, 581-592.\\

\noindent Rubin, D.B. Multiple Imputation for Nonresponse in
Surveys. {\it New York, Wiley}, 1987.\\

\noindent Rubin, D.B., Schenker, N. Multiple imputation for interval
estimation from simple random samples with ignorable nonresponse
{\it  J. Am. Statist. Assoc.}, 81, 1986, 366-374.\\

\noindent Schafer, J.L., Analysis of incomplete multivariate data.
London, {\it Chapman \& Hall}, 1997.

 \noindent SAS/STAT User's Guide, Vol. 2, Cary, NC, {\it SAS Institute Inc.}, 2002. \\

 \noindent Verbeke, G., Molenberghs, G. Linear mixed models for longitudinal
data. New York, {\it Springer}, 2000. \\

 \noindent Zhang, B. Confidence
intervals for a distribution function in the presence of auxiliary
information. {\it Computational statistics and data analysis}, v.21,
1996, pp 327-342.

\section*{\begin{center} Appendix: Large Sample Properties \end{center}}

If $\mathbf{K}_{ij}$, $\mathbf{K}^{*}_{ij}$ and $Cov \left ( \hat
\Theta_{ij}, \hat \Theta_{ij} \right)$ are known and there exists
$\left(\mathbf{K}^{*}_{ij}\right)^{-1}$, then the estimator (1) can
be calculated and the asymptotic properties of the estimator (1)
described by the following result.

 \textbf{Proposition 1.}  Let us consider the vectors $\hat\xi_{ij} \equiv
\sqrt{J_{ij}}\left(\hat\Theta_{ij}-\Theta_{ij}\right)$,
$\hat\psi_{ijs} \equiv \sqrt{J_{ij}}\left(\hat{\cal{B}}_{ijs} -
{\cal{B}}_{ijs}\right)$, $\hat\zeta_{ijs} \equiv \sqrt{J_{ijs}}
\left(\tilde{\cal{B}}_{ijs} - {\cal{B}}_{ijs}\right)$,
$s=1,...,S^{*}$ (for simplicity we omit ij-subscript in further
notation) with the following properties:

1)  $\hat\xi \longrightarrow \xi$,
 in distribution, as $J \rightarrow +\infty$.
 Also assume $E\left(\xi\right) \equiv \textbf{0}$ and all
elements composing covariance matrix
$Cov\left(\xi,\xi\right)\equiv\mathbf{C}^{\left(\xi,\xi\right)}$ are
finite.

2) $\hat\zeta^{(s)} \longrightarrow \zeta^{(s)}$, in distribution,
as $J_{s}\rightarrow +\infty$. Also assume
$E\left(\zeta^{(s)}\right)\equiv \mathbf{0}$ and all elements
composing covariance matrix $Cov\left(\zeta^{(s)},\zeta^{(s)}\right)
\equiv \mathbf{C}^{\left(\zeta,\zeta\right)}_{ss}$ are finite, for
all $s=1,...,S^{*}$.

3) $\hat\psi^{(s)} \longrightarrow \psi^{(s)}$, in distribution, as
$J \rightarrow +\infty$. Also assume $E(\psi^{(s)}) \equiv
\textbf{0}$ and all elements composing covariance matrices $Cov
\left(\psi^{(s)},\psi^{(q)} \right) \equiv
\mathbf{C}^{\left(\psi,\psi\right)}_{sq}$ and
$Cov\left(\xi,\psi^{(s)}\right)
\equiv \mathbf{C}^{\left(\xi,\psi\right)}_s$ are finite, for all $s,q=1,...,S^{*}$. \\

If $det\left(\mathbf{K}^{*}\right)>0$ and
$\frac{\sqrt{J}}{\sqrt{J_s}} {\longrightarrow} w_s \in [0,+\infty)$,
as $J$ and/or $J_s$ go to $+\infty$, then $\hat\eta \equiv \sqrt{J}
\left( \hat\Theta - \Theta \right) $ converges to a random vector
$\eta$ with $E\left(\eta\right)\equiv \mathbf{0}$ and
 $Cov\left(\eta,\eta\right) \equiv \mathbf{C}^{\left(\eta,\eta\right)}
 \equiv
\mathbf{C}^{\left(\xi,\xi\right)}-\mathbf{C}^{(\xi,\psi)} \left(
\mathbf{C}^{(\psi+\zeta)} \right)^{-1} \left(
\mathbf{C}^{(\xi,\psi)}\right)^{T}$, where matrices
$\mathbf{C}^{(\xi,\zeta)}$ and $\mathbf{C}^{(\psi+\zeta)}$ are
combined from the other matrices
$\mathbf{C}^{(\xi,\zeta)}=\|\mathbf{C}^{(\xi,\zeta)}_s\|_{s=1,...,S^{*}}$
and $\mathbf{C}^{(\psi+\zeta)}=\|\mathbf{C}^{(\psi,\psi)}_{sq}
+I_{[s=q]}w_s^2\mathbf{C}^{(\zeta,\zeta)}_{ss}\|_{s,q=1,...,S^{*}}$. \\

\textbf{Proof.} Taking into consideration that $\hat\Theta$ is an
unbiased estimator of $\Theta$ have $E\left(\eta\right) = \sqrt{J}
\left(E\hat\Theta - \Theta \right) = \mathbf{0}$. Hence,
$E\left(\eta\right) = \mathbf{0}$.

From (2) have $\mathbf{C}^{(\eta,\eta)} = J
\left(Cov\left(\hat\Theta,\hat\Theta\right)- \mathbf{K} \left(
\mathbf{K}^{*} \right)^{-1} \left( \mathbf{K} \right)^{T} \right)$.

Applying the facts

(1) $J Cov\left(\hat\Theta,\hat\Theta\right)$ converges to
$\mathbf{C}^{(\xi,\xi)}$, as $J$ goes to $+\infty$,

(2) $\mathbf{{K}} \left( \mathbf{K}^{*} \right)^{-1}$ converges to
$\mathbf{C}^{(\xi,\psi)}
\left(\mathbf{C}^{(\psi+\zeta)}\right)^{-1}$, as
$\frac{\sqrt{J}}{\sqrt{J_s}}$ goes to $w_s$, as $J$ and/or $J_s$ go
to $+\infty$, and

(3) $J\left(\mathbf{K}\right)^{T}$ goes to
$\left(\mathbf{C}^{(\xi,\psi)}\right)^T$, as $J$ goes to $+\infty$,

conclude $\mathbf{C}^{(\hat\eta,\hat\eta)}$ converges to
$\mathbf{C}^{(\eta,\eta)}$. \textbf{Q.E.D.} \\


In the expression for $\mathbf{C}^{(\eta,\eta)}$ the term
$\mathbf{C}^{(\xi,\psi)} \left(\mathbf{C}^{(\psi+\zeta)}\right)^{-1}
\left(\mathbf{C}^{(\xi,\psi)}\right)^{T}$ consists of quadratic
forms and corresponds to the decrease of the original dispersion
ellipsoid. Applying different quadratic forms (defined by risk
function) to $\mathbf{C}^{(\eta,\eta)}$ the term
$\mathbf{C}^{(\xi,\psi)} \left( \mathbf{C}^{(\psi+\zeta)}
\right)^{-1} \left( \mathbf{C}^{(\xi,\psi)}\right)^T$ defines
different non-negative numbers showing asymptotic improvement of
used risk function.

In the Proposition 1 the cases when there exists $w_s=+\infty$ were
not considered because as only $w_s=+\infty$ information from
$s^{th}$-subsample on $(i+1)$ level is overwhelmed by information in
$ij$-subsample and cannot improve the asymptotic properties of the
estimators derived from $ij$-subsample. In this case
$s^{th}$-subsample on $(i+1)$ level should be excluded from
consideration.

Another extreme situation appears when $w_s$ is equal to 0 which
corresponds to incorporating information of exact knowledge. In the
case $\tilde{\cal{B}}$ is known with zero variance.

Proposition 1 defines the asymptotic properties of the estimator (1)
but this estimator cannot be used in a number of practical cases
because $\mathbf{{K}} \left( \mathbf{K}^{*} \right)^{-1}$ usually is
not known. In this case the estimator $\hat{\tilde{\Theta}}$,
obtained in (4) by substitution $\mathbf{{K}} \left( \mathbf{K}^{*}
\right)^{-1}$ on $\mathbf{\hat{K}} \left( \mathbf{\hat{K}}^{*}
\right)^{-1}$, should be used. The asymptotic properties of
$\hat{\tilde{\Theta}}$ is described as follows.
\\

\textbf{Proposition 2.} Suppose the assumptions of Proposition 1
hold and every element of $J
\left(\widehat{Cov}\left(\hat\Theta,\hat\Theta\right) -
Cov\left(\hat\Theta,\hat\Theta\right)\right)$, $J
\left(\hat{\mathbf{K}} - \mathbf{K}\right)$, and  $J
\left(\hat{\mathbf{K}}^{*} - {\mathbf{K}}^{*} \right) $ converges to
some random variable with mean zero and finite variance.

Then $\sqrt{J} \left(\hat{\tilde{\Theta}} - \Theta \right)$
converges in distribution to $\eta$, as $J \longrightarrow +\infty$,
where $\eta$ defined in Proposition 1.

\textbf{Proof}.

Notice that $\sqrt{J} \left(\hat{\tilde{\Theta}} - \Theta \right)$
differs from $\sqrt{J} \left(\tilde{\Theta} - \Theta \right)$ only
by applying $\mathbf{\hat K}$ and $\mathbf{\hat {K}^{*}}$ instead of
$\mathbf{K}$ and $\mathbf{{K}^{*}}$.

From the fact that linear combinations of elements of $\mathbf{\hat
K}$ and $\left(\mathbf{\hat K}^{*}\right)^{-1}$ are continuous
functions and all these elements converge in probability to their
true values on the basis of Theorem 5.5.4 (Casella and Berger, 2002,
p. 233) conclude $\mathbf{\hat K}\left(\mathbf{\hat
K}^{*}\right)^{-1}$ converges in probability to
$\mathbf{K}\left(\mathbf{K}^{*}\right)^{-1}$.

Now from Slutsky's Theorem (Casella and Berger, 2002, p. 239)
conclude $\sqrt{J} \left(\hat{\tilde{\Theta}} - \Theta \right)$
converges in distribution to $\eta$.
Q.E.D. \\

\emph{Remark:} For the cases when all $w_s=0$ estimator (4) becomes
the same as the estimator derived by method of correlated processes
(Pugachev,1973) and has the same asymptotical properties as the
empirical likelihood estimator in the presence of auxiliary
information (Zhang,1996).

\end{document}